\documentclass[11pt]{article}
\pdfoutput=1

\usepackage{amsmath,amsfonts}
\usepackage[lmargin=1.5in,rmargin=1.5in,top=1in,bottom=1in]{geometry}
\usepackage[tiny,compact]{titlesec}
\usepackage{url}
\usepackage{graphics}
\usepackage{titlesec}
\usepackage{graphicx}
\usepackage{textcase}
\usepackage{epstopdf}
\usepackage{url}
\usepackage{booktabs}
\usepackage{tabularx}
\usepackage{enumitem}
\usepackage{multicol}
\usepackage{sidecap}
\usepackage[font={scriptsize,it}]{caption}
\usepackage{color}   
\usepackage{hyperref}
\hypersetup{
    colorlinks=true, 
    linktoc=all,     
    linkcolor=blue,  
}
\usepackage{subfig}
\usepackage{bibentry}
\usepackage{framed}

\graphicspath{{./}{./figs/}}

  \def\clap#1{\hbox to 0pt{\hss#1\hss}}

\providecommand{\mat}[1]{\boldsymbol{#1}}%
\renewcommand{\vec}[1]{\mathbf{#1}}

\providecommand{\mX}{\ensuremath{\mat{X}}}

\providecommand{\va}{\ensuremath{\vec{a}}}

\providecommand{\vf}{\ensuremath{\vec{f}}}

\providecommand{\vu}{\ensuremath{\vec{u}}}

\providecommand{\vw}{\ensuremath{\vec{w}}}
\providecommand{\vx}{\ensuremath{\vec{x}}}

\newcommand{\bmat}[1]{\begin{bmatrix}#1\end{bmatrix}}

\newcommand{\argmin}[1]{\underset{#1}{\operatorname{argmin}}}

\graphicspath{{./}{./figs/}}
\usepackage{epstopdf}

\titleformat*{\section}{\Large\bfseries}
\titleformat*{\subsection}{\large\bfseries}
\titleformat*{\subsubsection}{\itshape}

\newcommand{\bi}{\begin{itemize}}
\newcommand{\ei}{\end{itemize}}

\begin{document}

\setcounter{page}{1}
\pagenumbering{arabic}

\begin{center}
\Large{\textsc{A Quick-and-Dirty Check\\ for a One-dimensional Active Subspace}}\\[5mm]
Paul G.~Constantine\\
Colorado School of Mines, Golden, CO\\[5mm]
\url{paul.constantine@mines.edu}\\
\url{http://inside.mines.edu/~pconstan}
\end{center}

\begin{abstract}
Most engineering models contain several parameters, and the map from input parameters to model output can be viewed as a multivariate function. An active subspace is a low-dimensional subspace of the space of inputs that explains the majority of variability in the function. Here we describe a quick check for a dominant one-dimensional active subspace based on computing a linear approximation of the function. The visualization tool presented here is closely related to \emph{regression graphics}~\cite{cook2009introduction,cook2009regression}, though we avoid the statistical interpretation of the model. This document will be part of a larger review paper on active subspace methods. 
\end{abstract}

We can run a quick-and-dirty check to see if a given model has a dominant one-dimensional active subspace---i.e., a direction in the space of model inputs that explains the majority of the variability in the model. After we spell out the steps of this check, we'll describe what it's looking for. Here's what we need to get started:
\begin{itemize}
\item {\it A scalar quantity of interest}. The model should have a single output $f=f(\vx)$ that depends on the $m$ input parameters denoted $\vx=(x_1,\dots,x_m)$. For example, this might be the lift of an airfoil, the maximum power of a battery, or the average temperature over a surface.
\item {\it A range for each input parameter}. Find a lower bound $x^l_i$ and an upper bound $x^u_i$ such that $x^l_i\leq x_i\leq x^u_i$ for each input parameter $x_i$. This may be a 10\% perturbation around some nominal value. 
\item {\it Enough time and computing power to run a few model evaluations}. To fit the linear model below, we need to evaluate the model $f$ about three-to-four-$m$ times. 
\end{itemize}
Given these, we run the following procedure.
\vspace{1ex}

\begin{figure}[htbps]
\subfloat[Scatter plot]{
\includegraphics[width=0.4\textwidth]{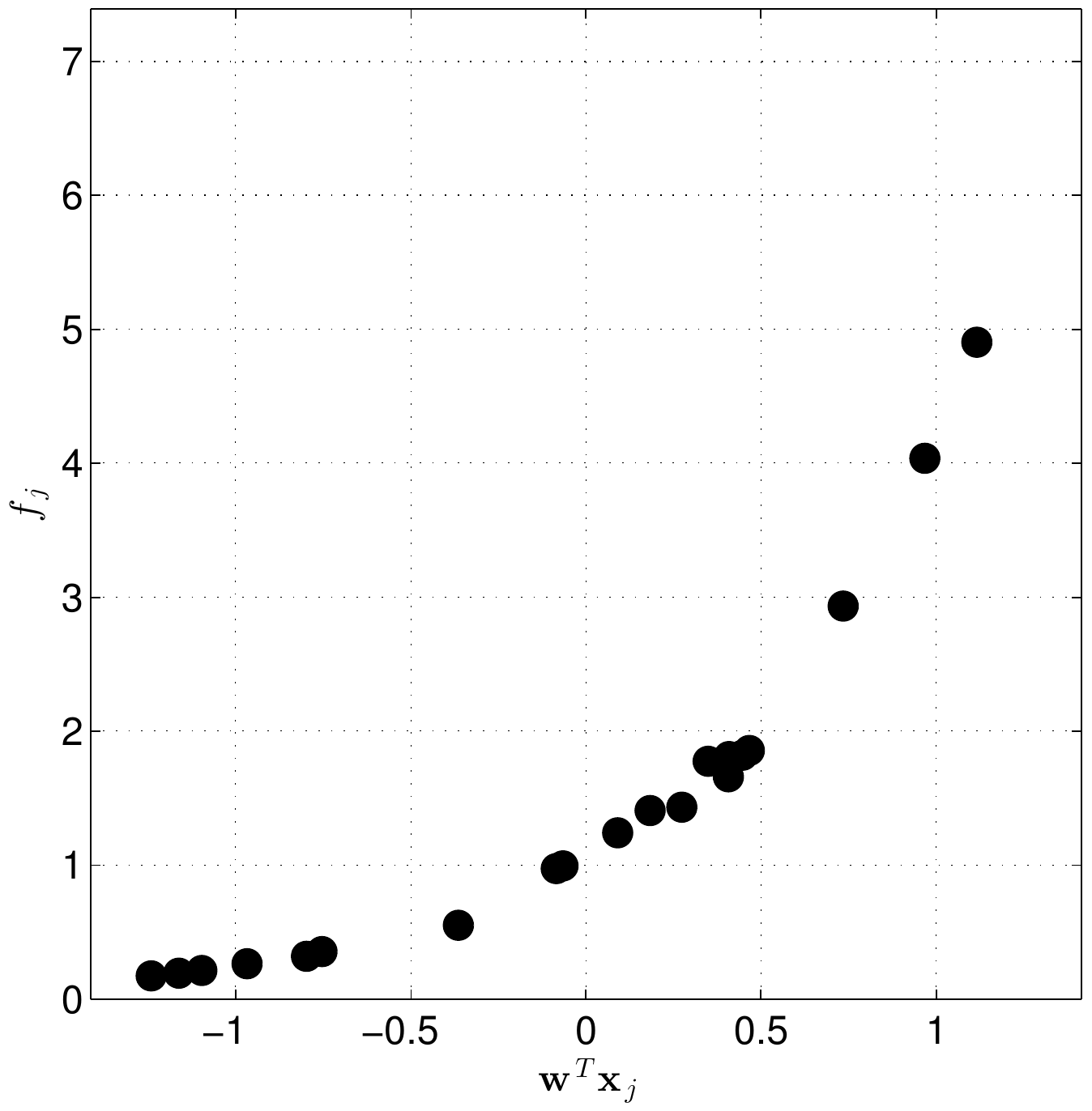} 
\label{fig:scatter}
}\;
\subfloat[Rotated surface plot]{
\includegraphics[width=0.5\textwidth]{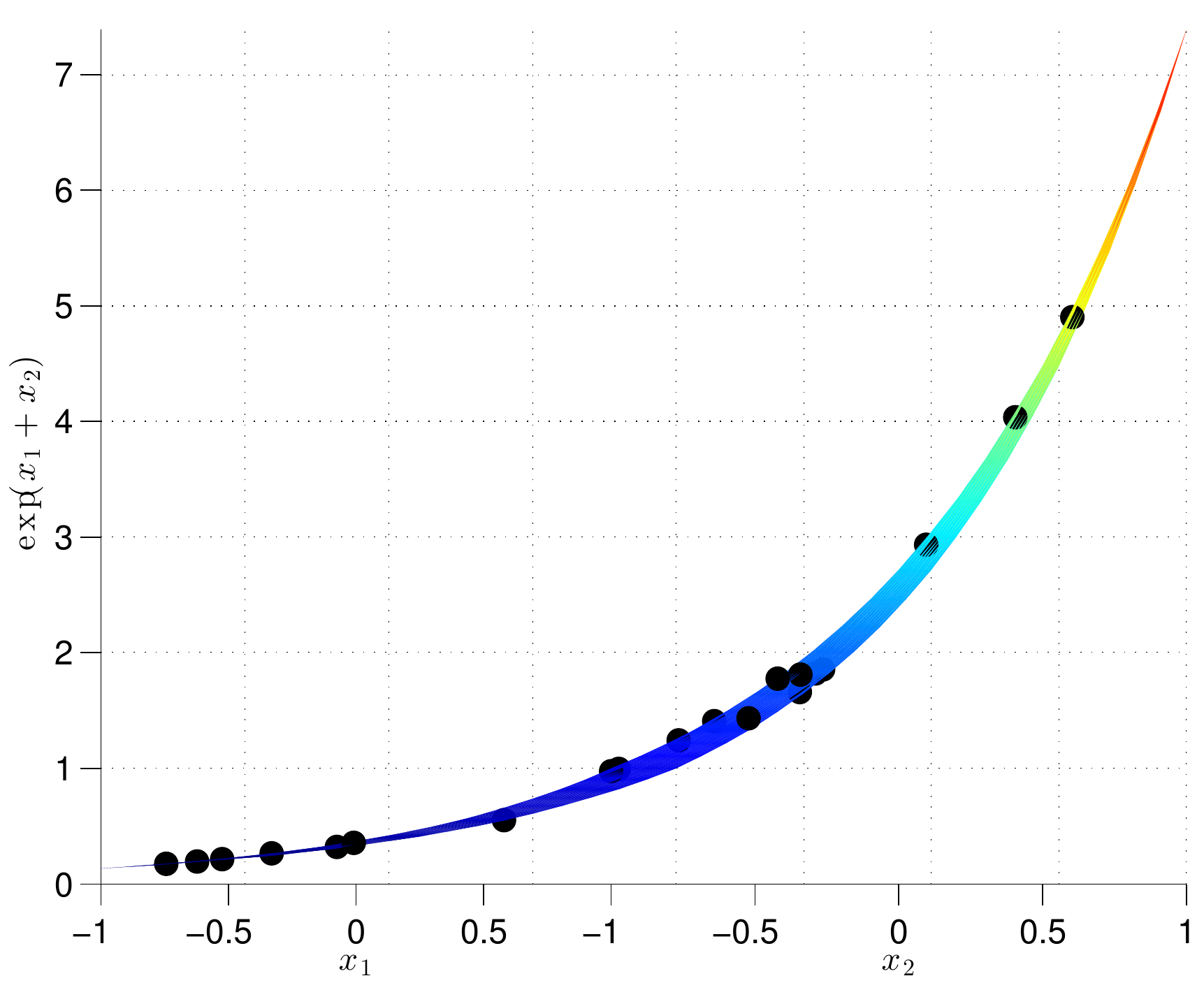}
\label{fig:surf}
}
\caption{The leftmost figure is the scatter plot produced by the quick check for the active subspace applied to the function $\exp(x_1+x_2)$. The rightmost figure is a rotated surface plot of the same function.}
\label{fig:exp}
\end{figure}

\noindent{\it A quick-and-dirty check for an active subspace}
\vspace{-2mm}
\begin{enumerate}
\item Draw $N$ samples $\hat{\vx}_j$ uniformly from $[-1,1]^m$.
\item Let $\vx_j=\frac{1}{2}\left[(\vx^u-\vx^l)\cdot\hat{\vx}_j+(\vx^u+\vx^l)\right]$, where $\vx^l$ is the vector of lower bounds and $\vx^u$ is the vector of upper bounds; the dot operation $(\cdot)$ is component-wise multiplication.
\item Compute $f_j = f(\vx_j)$. 
\item Compute the coefficients of the linear regression model
\[
f(\vx) \;\approx\; \hat{a}_0 + \hat{a}_1\hat{x}_1 + \cdots + \hat{a}_m\hat{x}_m
\]
with least-squares
\[
\hat{\va}\;=\;\argmin{\vu}\; \|\hat{\mX}\vu - \vf\|_2^2,
\]
where 
\[
\hat{\mX} = \bmat{1 & \hat{\vx}_1 \\ \vdots & \vdots \\ 1 & \hat{\vx}_N},\quad
\hat{\va} = \bmat{\hat{a}_0 \\ \vdots \\ \hat{a}_m},\quad
\vf = \bmat{f_1\\ \vdots \\ f_N}.
\]
\item Let $\vw = \va'/\|\va'\|$ where $\va'=[\hat{a}_1,\dots,\hat{a}_m]^T$ is the last $m$ coefficients (i.e., the gradient) of the linear regression approximation.
\item Produce a scatter plot of $\vw^T\hat{\vx}_j$ versus $f_j$. 
\end{enumerate}
Here is a MATLAB implementation of this procedure, assuming there is a function \texttt{mymodel.m} that computes the model output $f$ given inputs and vectors \texttt{xl} and \texttt{xu} containing the lower and upper bounds of the input parameters.
\begin{verbatim}
% Produce a scatter plot to check for a 
% dominant 1-d active subspace.
N = 4*m;
Xhat = 2*rand(m,N);
X = 0.5*(repmat(xu-xl,1,N).*Xhat + repmat(xu+xl,1,N));
for j=1:N, f(j) = mymodel(X(:,j)); end
ahat = [ones(N,1) Xhat'] \ f';
w = ahat(2:m+1)/norm(ahat(2:m+1));
plot(Xhat*w,f,'o');
\end{verbatim}
To demonstrate, we apply this procedure to $f(\vx)=\exp(x_1+x_2)$ defined on the square $[-1,1]^2$. Figure \ref{fig:scatter} shows the scatter plot using $N=20$. Notice the apparent relationship between evaluations $\vw^T\vx_j$ and $f_j$. In fact, this particular $f$ varies entirely along a one-dimensional subspace defined by the vector $[1, 1]^T$. 

The scatter plot in Figure \ref{fig:scatter} is equivalent to sampling the three-dimensional surface plot---where the $z$ coordinate is $\exp(x_1+x_2)$---and rotating the plot such that the evaluations $f_j$ appear to depend on one variable instead of two. The rotated surface plot along with the $f_j$ is shown in Figure \ref{fig:surf}. This interpretation is valid more generally for functions of several variables. The check for an active subspace finds an angle from which to view the data set $\{\vx_j,f_j\}$ such that a trend emerges in the $f_j$ (if such a trend exists) as a function of a new variable $y=\vw^T\vx$.

If, instead of the linear regression procedure, we choose the vector $\vw$ as the $i$th column of the $m\times m$ identity matrix, then we create a scatter plot of $f_j$ versus the $i$th component of $\vx_j$. Such scatter plots are described in Section 1.2.3 of \cite{saltelli2008global}. However, the vector produced by the regression is special. It is the normalized gradient of the linear regression approximation to $f(\vx)$. If a global, monotonic trend is present in the $f$---and if the number of samples $N$ is sufficiently large---then the vector $\vw$ reveals the direction in the parameter space of this trend. Additionally, the components of $\vw$ measure the relative importance of the original parameters $\vx$. For the function $\exp(x_1+x_2)$, we expect the two components of $\vw$ to be roughly the same. 

This procedure is a subjective diagnostic test based on visualization. The viewer must use the scatter plot to judge if such a trend is present. A trend may be obvious, subtle, or absent. For complex functions, it is unlikely that one can know if such a trend exists without running the check. Fortunately, the check is pretty cheap! 

Once a trend is identified, there are several ways to exploit it. For example, if one wants to maximize or minimize $f(\vx)$, then the vector $\vw$ is a direction in the parameter space along which we expect $f$ will increase or decrease. In particular, when the domain of $\vx$ is an $m$-dimensional hypercube, then the vector $\vw$ can identify one of the $2^m$ corners where $f$ is likely to be the largest. If $m$ is large and $f$ is expensive, then this strategy is preferred to evaluating $f$ $2^m$ times. 

If one is attempting to invert $f$ to find parameters that match a given measurement, then this trend can help identify sets of parameters that map to the measurement. If one wishes to predict $f$ at some other values of $\vx$, then this trend can help build a coarse surrogate model on a low-dimensional subspace of the parameter space. The trend may also be used to approximate an integral of $f$ over the parameter space. 

\noindent\textbf{Examples.}

\begin{figure}[htbps]
\subfloat[]{
\includegraphics[width=0.53\linewidth]{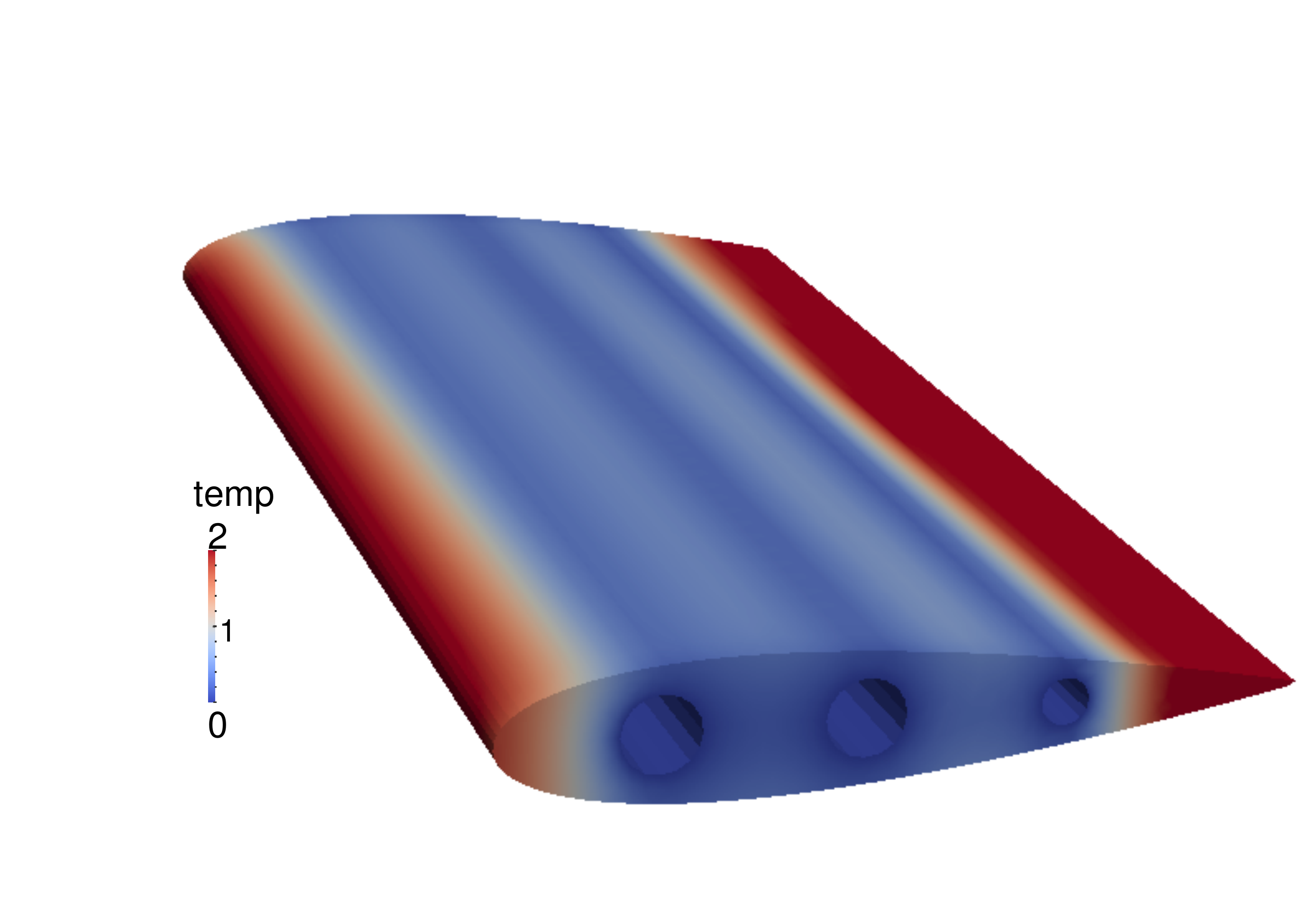}
}
\subfloat[]{
\includegraphics[width=0.37\linewidth]{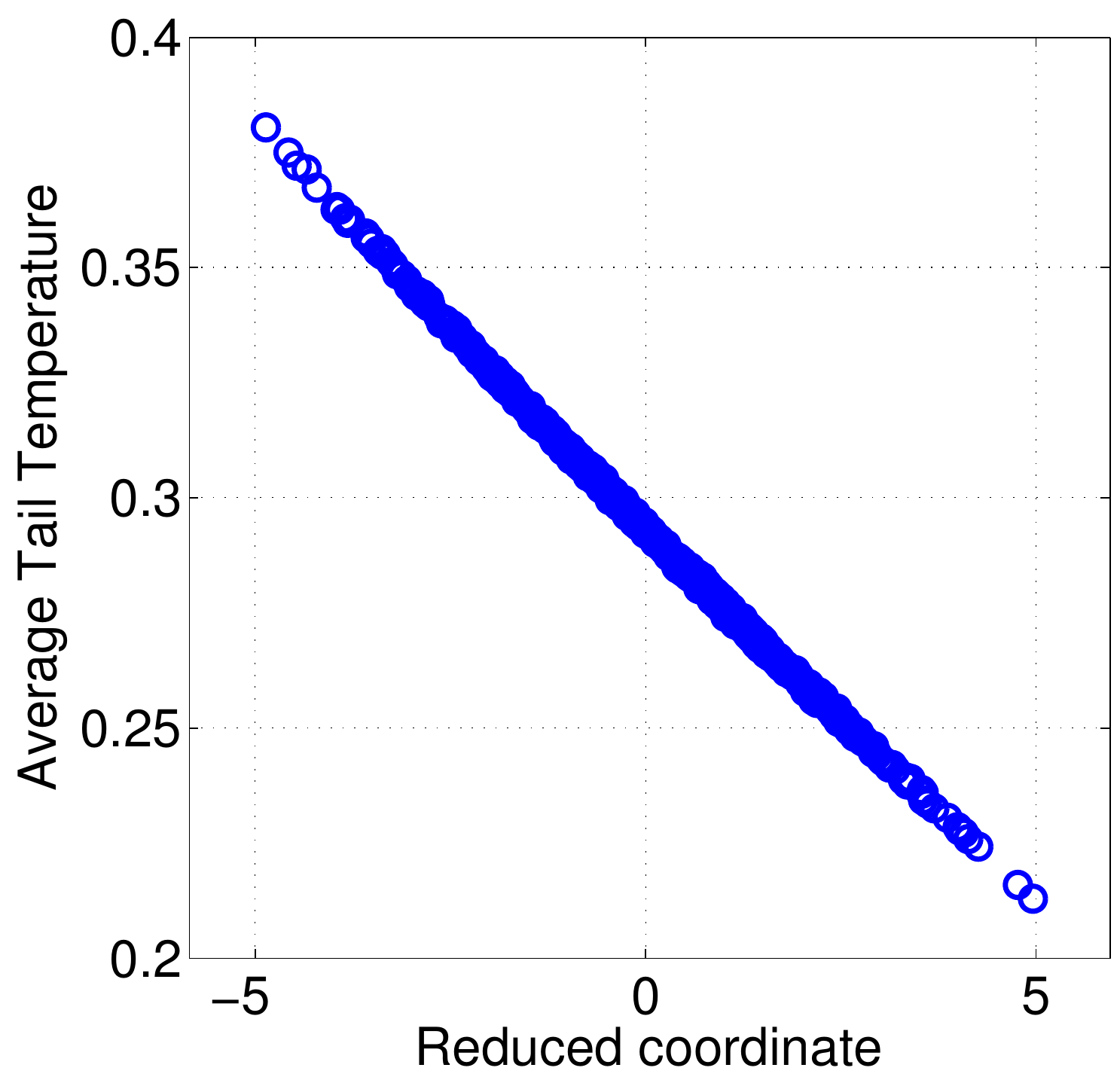}
}
\caption{We examined a model for heat transfer in a turbine blade given a parameterized model for the heat flux boundary condition representing unknown transition to turbulence. There are 250 parameters characterizing a Karhunen-Loeve model of the heat flux, and the quantity of interest is the average temperature over the trailing edge of the blade. The leftmost figure shows the domain and a representative temperature distribution. The rightmost figure plots 750 samples of the quantity of interest against the projected coordinate. The strong appearance of the linear relationship verifies the quality of the subspace approximation.}
\label{fig:blade}
\end{figure}

\begin{figure}[htbps]
\centering
\subfloat[]{
\includegraphics[width=0.5\linewidth]{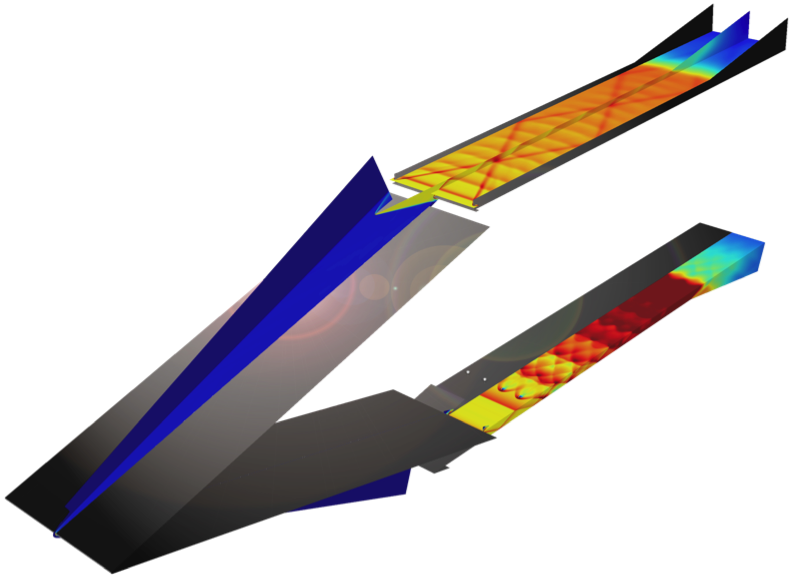}
}
\subfloat[]{
\includegraphics[width=0.4\linewidth]{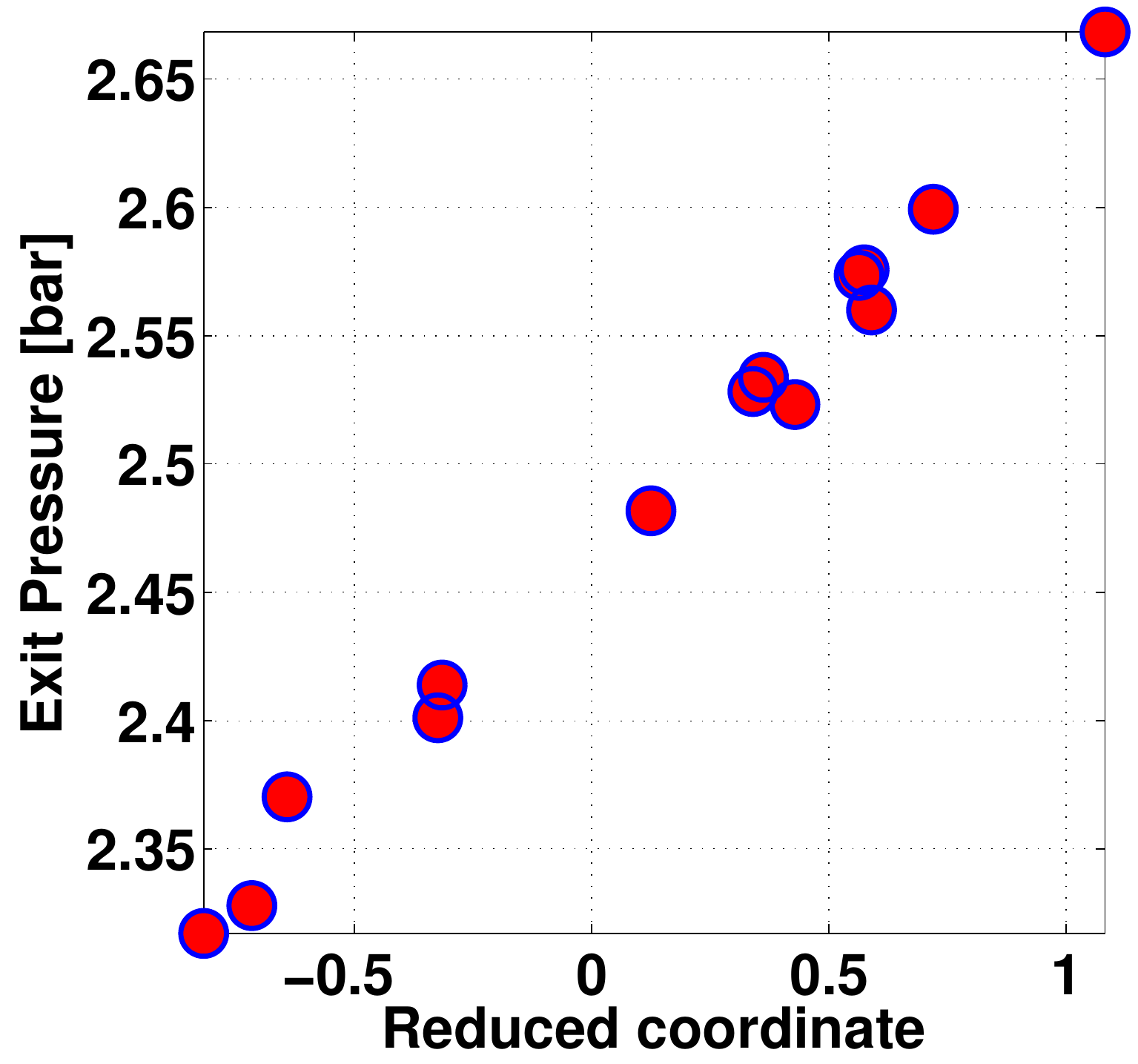}
}
\caption{We applied the active subspace check to quantify the margins of safety in a multiphysics model of a hypersonic scramjet vehicle. There are seven parameters characterizing the inflow boundary conditions of a reacting, compressible channel flow. The quantity of interest was the exit pressure at the exit of the combustor. The left figure shows a representative compressible flow computation. The right figure shows the scatter plot of 14 samples of the exit pressure.}
\label{fig:hyshot}
\end{figure}

\begin{figure}[htbps]
\centering
\subfloat[]{
\includegraphics[width=0.45\linewidth]{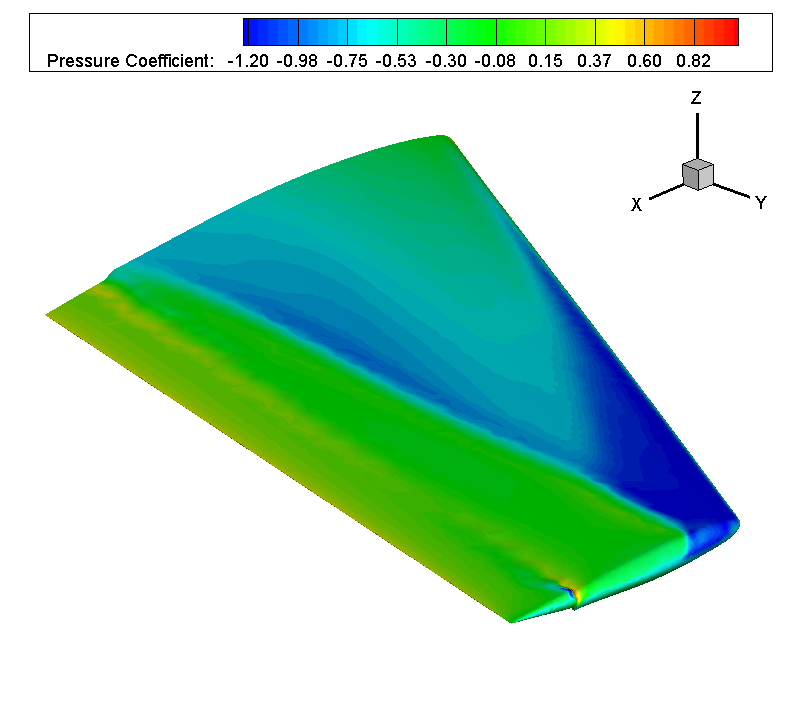}
}
\subfloat[]{
\includegraphics[width=0.45\linewidth]{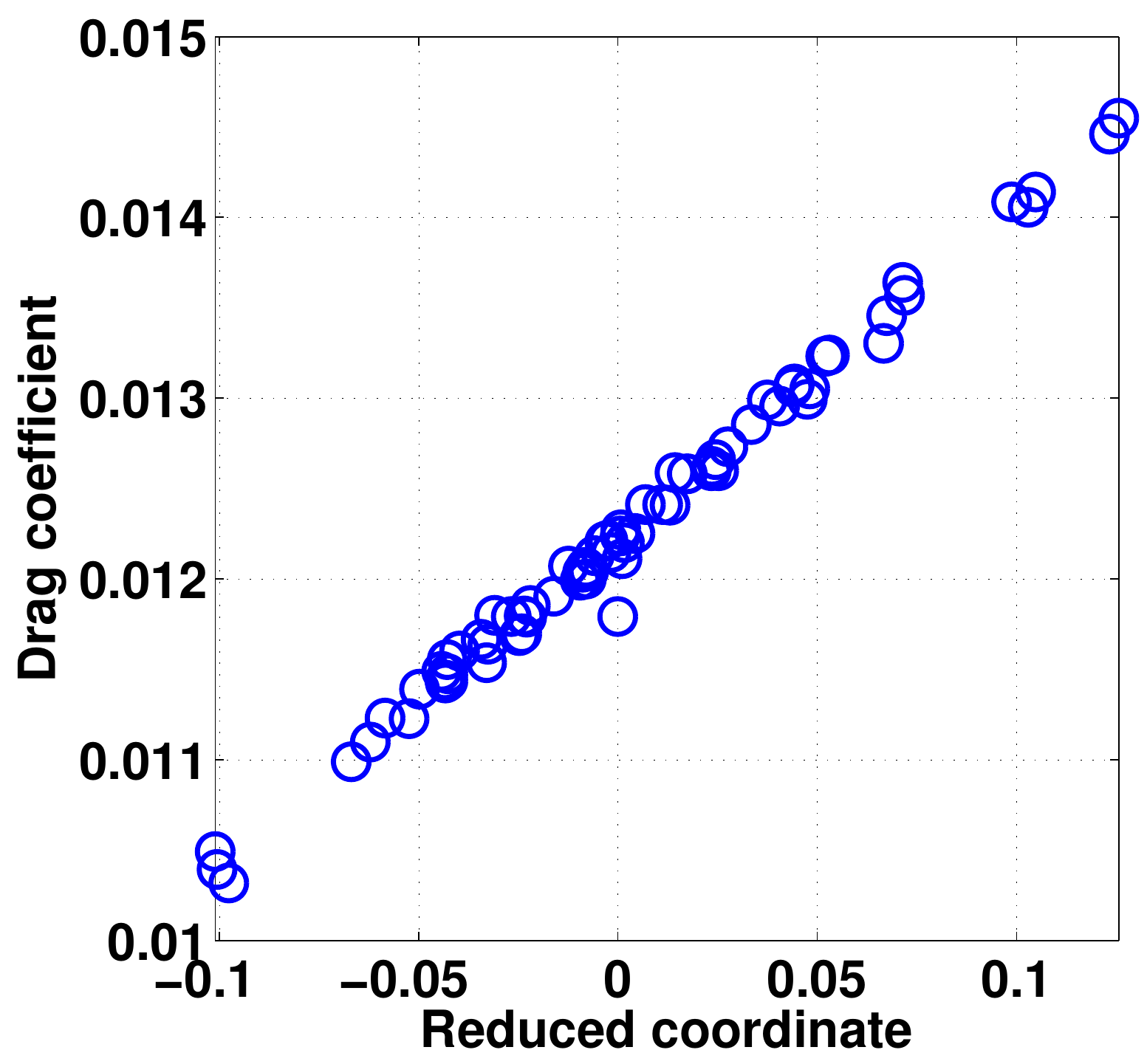}
}
\caption{We applied the active subspace check to a design optimization problem in transonic wing design. There were 60 design variables, and the quantity of interest is the drag coefficient computed via computational fluid dynamics. The left figure shows the a representative pressure field for a given design. The right figure shows 70 evaluations of the drag coefficient plotted against the reduced coordinate of the one-dimensional active subspace from the 60-dimensional design space. The nearly linear relationship implies that optimization of the drag coefficient is much easier than initially thought.}
\label{fig:oneram}
\end{figure}

The check for the active subspace has worked surprisingly well for several complex systems. In particular, it has identified a direction in the parameter space that corresponds to a global monotonic trend in the output. Three representative systems and the corresponding plots are shown in Figures \ref{fig:blade}, \ref{fig:oneram}, and \ref{fig:hyshot} with a brief description in the captions. More details on these systems can be found in \cite{constantine2013active}, \cite{Lukaczyk2014}, and \cite{icossar2013}, respectively. We are actively seeking more examples. 

There are two common features of these complex models that may explain the perceived trends:
\begin{enumerate}
\item When a scientist informally describes the effect of a parameter on a model, it often sounds like the output is monotonic with respect to changes in the parameter; more of $x_3$ leads to less $f$. This may suggest a global monotonic trend in the model as a function of the parameter.
\item A model is often given a nominal parameter value corresponding to a theoretical or measured physical case. Back-of-the-envelope estimates of the uncertainty in the input parameters often take the form of a 5-10\% perturbation about the nominal value. Despite the complexity of the system, this type of perturbation often results in little change in the model. More importantly, it results in small rate of change. In such cases, the linear model represents the global trends well within the range of the input perturbations.  
\end{enumerate}
These explanations are based on experience with several models. Of course, there are many complex systems that do not behave this way. In such cases, the check for the active subspace may lead to nothing useful. However, given the low cost of the check and the potential insights into the model, we think it is well worth the effort.

\bibliographystyle{siam}
\bibliography{asm-primer}

\end{document}